\documentclass{article}
\usepackage{graphicx}
\usepackage{amsmath}
\usepackage{amsthm}
\usepackage{amsfonts}
\usepackage{amssymb}
\usepackage{enumerate}
\usepackage{tikz}
\usepackage{youngtab}
\usepackage{mathtools}
\usepackage{blkarray}
\usepackage{empheq}
\usepackage{subfig}
\usepackage{floatrow}
\usepackage{siunitx}
\usepackage{booktabs}
\DeclarePairedDelimiter\ceil{\lceil}{\rceil}

\DeclarePairedDelimiter\floor{\lfloor}{\rfloor}

\newcommand*{\horzbar}{\rule[.5ex]{2.5ex}{0.5pt}}
\usepackage{pgfplots}
\usepackage{tabularx}

\usepackage{hyperref}
\usepackage{nicematrix}

\DeclareSymbolFont{largesymbolsA}{U}{jkpexa}{m}{n}
\SetSymbolFont{largesymbolsA}{bold}{U}{jkpexa}{bx}{n}
\DeclareMathSymbol{\varprod}{\mathop}{largesymbolsA}{16}

\usepackage{caption}
\usepackage{algorithm}
\usepackage{algorithmic}

\usepackage[english]{babel}
\usepackage[square,numbers]{natbib}
\title{Fast Symbolic Integer-Linear Spectra}
\author{Jonny Luntzel, Abraham Miller}
\date{}

\pgfplotsset{compat=1.18}
\begin{document}

\maketitle

\section*{Abstract}
Here we contribute a fast symbolic eigenvalue solver for matrices whose eigenvalues are $\mathbb{Z}$-linear combinations of their entries, alongside efficient general and stochastic $M^{X}$ generation. Users can interact with a few degrees of freedom to create linear operators, making high-dimensional symbolic analysis feasible for when numerical analyses are insufficient.

\section*{Introduction}
$\mathbb{Z}$-linear eigenvalue matrices are a class of matrices which have monomial entries, and eigenvalues expressed as sums of integer multiples of their entries. Kenyon\cite{ZLE} and colleagues explore a subset of such matrices\cite{OS}, and describe a construction paraphrased below:

\begin{itemize}
    \item Select a finite strict partial order $X$ (fig. \ref{fig:o}) 
    \item Determine size-n permutations (linear extensions) $\phi$ which satisfy $X$ 
    \item Create a matrix whose (i,j) entries are $\phi_{j}\phi_{i}^{-1}(id)$
    \item Apply the ascent-descent function $\epsilon$ on each entry (e.g. $132\rightarrow10$)
\end{itemize}

\begin{figure}[H]
\begin{center}
    $X=(2<1, 2<3) \;\;\rightarrow \;\;\phi=\{213,312\}$
\end{center}
\begin{center}
\scalebox{0.62}{
    \begin{tikzpicture}
    \node[shape=circle,draw=black] (A) at (0,0) {1};
    \node[shape=circle,draw=black] (C) at (1,2) {3};
    \node[shape=circle,draw=black] (B) at (2,0) {2};

    \path[->,draw=black]
    (B) edge node {} (C)
    (B) edge node {} (A)
    ;
    \end{tikzpicture}
    }
\[
\begin{pNiceMatrix}[first-col,first-row]
  & 213 & 312 \\
213 & 123 & 321\\
231 & 321 & 123 
\end{pNiceMatrix} \rightarrow 
\begin{pNiceMatrix}
11 & 00\\
00 & 11
\end{pNiceMatrix}\rightarrow 
\begin{pNiceMatrix}
a_{1} & a_{2}\\
a_{2} & a_{1}
\end{pNiceMatrix}
\]
\end{center}
\caption{A partial order and its corresponding one-line permutations. Entries of $M^X$ are row-inverses composed with column permutations, reduced to $01$-strings by $\epsilon$. Strings act as monomials in the resultant symbolic matrix.}
\label{fig:o}
\end{figure}

\section*{Matrix Generation}
Users can generate arbitrary $M^{X}$ with partial order inputs\footnote{\url{https://github.com/orgs/symeig/repositories}, 2024.\cite{lapack}\cite{blas}\cite{sympy}\cite{mpmath}\cite{MPLAPACK}}, owing to the \hyperref[sec:gen]{script's}  dynamic programming approach which generates legal permutations from ground-up instead of filtering from the set of possible permutations.
\newline\\
Furthermore, a subset of partial orders always yield stochastic matrices with predictable properties. In order to guarantee stochasticity we \hyperref[sec:fgen]{restrict} to partial orders which are fixed disjoint blocks of local transpositions, and directly take products on the possible binary sequences rather than from permutations themselves to assign monomial entries, described below:

\begin{itemize}
    \item Select a partial order $X$ (determined by disjoint block lengths)
    \item Generate the corresponding $\epsilon$-filtration of possible $\pm$ sequences ($\ll 2^{n-1}$)
    \item Create a matrix whose $(i,j)$ entries are $\epsilon_{\phi_{j}\phi_{i}}$ computed directly from the filtration
\end{itemize}
Rather than a partial order, Users specify the dimension $n$ with a Fibonacci-factorization of $n$ in order to specify the contribution of each local block to the total dimension. Factors must be Fibonacci because they indirectly specify the length of local swap chains added to the matrix's partial order, and such $k$-length chains yield $a_k = a_{k-1} + a_{k-2}$ permutations. For example, $\{13,2\} \rightarrow n=26$.

\begin{figure}[H]
\begin{tikzpicture}
    \node[] (w) at (-0.5,0) {$\{$};
     \node[shape=circle,draw=black] (A) at (0,0) {1};
    \node[shape=circle,draw=black] (B) at (1.5,0) {2};
    \node[shape=circle,draw=black] (C) at (3,0) {3};
    \node[shape=circle,draw=black] (D) at (4.5,0) {4};
    \node[shape=circle,draw=black] (E) at (6,0) {5};
    \node[] (x) at (6.5,0) {$\}$};
    
     \node[] (x) at (7.5,0) {$\{$};
    \node[shape=circle,draw=black] (F) at (8,0) {6};
    \node[shape=circle,draw=black] (G) at (9.5,0) {7};
     \node[] (x) at (10,0) {$\}$};
5, 6, 5, 7, 4, 6, 4, 7
    \path[->,draw=black]
    (C) edge [bend right] node[below] {} (A)
    (D) edge [bend left] node[below] {} (A)
    (D) edge [bend right] node[below] {} (B)
    (E) edge [bend right] node[below] {} (C)
    (E) edge [bend left] node[below] {} (B)
    (F) edge [bend left] node[below] {} (E)
    (G) edge [bend right] node[below] {} (E)
    (F) edge [bend right] node[below] {} (D)
    (G) edge [bend left] node[below] {} (D);

\end{tikzpicture}
\caption{Two disjoint local transposition chains where $n=16$ and $\{8,2\}$ dimension factors. Chain lengths of length $k$ generate $fib(k)$ permutations, and matrix dimension $n$ is a product of the number of permutations each chain produces.}
\label{fig:dbo}
\end{figure}

\noindent Partial orders are closed under concatenation of local chains, so while this construction is a strict subset, it guarantees the property holds. $M^{X}$ with disjoint block structure guarantees $\mathbb{Z}$-linear eigenvalues and provides users with a simple and well-leveraged interface, while pre-computing the direct filtration increases speed \& scalability. 

\begin{figure}[H]
\includegraphics[width=.56\textwidth]{0702.jpg}
\caption{Projective barycentric coordinate for a conservative $n=3$ system where $a+b+c=1$. Allows negative masses and extension to infinity (points marked with $\pm$ are $\pm \infty$ in a given variable for $abc$).}
\label{fig:oo}
\end{figure}

\noindent Users provide a small number of degrees of freedom in a highly flexible subset of $M^X$. Such matrices are desirable for Markov-type systems which are non-dissipative but not necessarily reversible (fig. \ref{fig:oo}). The general matrix generation (gen.py) lets users work directly with the relationships between the partial order and its corresponding matrix, while the disjoint, local transposition matrix generation (fgen.py) takes the matrix dimension as input and generates matrices with a chain structure, removing the interface to partial orders and ensuring stochasticity.\newline\
\begin{center}
Encoding Scheme
\end{center}
\begin{figure}[H]
\begin{enumerate}
  \item $x_{i} = 10^{3(k-i)}$
  \item $\lambda_{t_{0}} = \underbrace{\color{black}00\color{red}0}_{\color{red}z}\color{black}|\underbrace{\color{black}00\color{orange}0}_{\color{orange}y}\color{black}|\underbrace{\color{black}00\color{violet}0}_{\color{violet}x}$, $\quad res = 444444444+\lambda(A_{x})$
  \item $res = \hdots\color{violet}\underbrace{\color{black}367}_{} \color{black} \rightarrow (7-4)*10^{0} + (6-4)*10^{1} + (3-4)*10^{2} = -77\color{violet}x$
\end{enumerate}
\caption{The eigenvalue solver's method for substituting symbolic calculations with numeric ones over $A_x$.}
\label{fig:encoding}
\end{figure}
(1) Associate monomials with a base power.\newline (2) Give each monomial a window of digits to accumulate, centered around a midpoint ($4$ here). \newline (3) Count each $x_i$ ($x,y,z$), output $c_{1}x + c_{2}y + c_{3}z$.
\newpage

\subsection*{Eigenvalue Computation}
\noindent Time complexity of numerical eigenvalue calculation with QR is $O(n^3)$ for dense matrices\cite{Strang}. Big O time complexity for symbolic eigenvalue calculation is extremely high and far less clear, suffering from complexities of symbolic root finding and intermediate expression swell \cite{ES}. Our algorithm achieves similar time complexity to the numerical case made possible by casting the symbolic problem into a numerical form\ref{fig:encoding}. 
\newline\\
Symbols used in elements of the matrix are encoded as unique imaginary power terms, resulting in pseudo-symbolic numerical values for each element, reminiscent of Gödel numbering, but here expressly with the purpose of computational performance. By exploiting the structure of $\mathbb{Z}$-linearity we ensure power terms do not mix\footnote{Symbolic mixing can occur from digit spillover. We mitigate this by centering around a midpoint and staggering power terms.}.
\newline\\
The time complexity is $O(n^3 \cdot b \log b)$ for dense matrices of dimension $n$ with $b$ bits of precision, where $b$ is in in the worst case a constant multiple of $n$. As a parallel algorithm however, the parallel span is only $O(n^3)$, since the impact of precision becomes constant as the number of processors $p \rightarrow n$ in the batched computation.
\newline\\
\noindent Each batch requires digit precision only equal to its batch size and can be run in parallel, masking out all but the smaller set of batched symbols in the input matrix before running the numeric solver. The respective partial eigenvalues are summed together after to achieve the final result.
\newline\\
We created implementations of the $\mathbb{Z}$LE eigenvalue algorithm in Python, Mathematica, and C++ to allow a range of accessibility and performance. Algorithm details are in the repository and \hyperlink{alg}{appendix}.\newline

  \begin{minipage}{0.7\textwidth}
  \centering
  Runtime Comparison\newline
    \includegraphics[width=\linewidth]{runtime_gh.jpg}
  \end{minipage}
  \hfill
  \begin{minipage}{0.27\textwidth}
    \centering
    C++ (s)\newline
    \scalebox{0.8}{
        \begin{tabular}{ c c c}
          n & $\mu$ & $\sigma$\\
          13 & 0.010 & 0.047e-2\\
          21 & 0.027 & 0.047e-2 \\
          34 & 0.091 & 0.082e-2 \\
          55 & 0.319 & 0.002 \\
          89 & 1.723 & 0.057 \\
          144 & 6.381 & 0.246 \\
          233 & 43.810 & 0.480 \\
          377 & 158.374 & 1.590 \\
          610 & 690.556 & 25.591 \\
          987 & 3007.097 & 108.985
        \end{tabular}    
    }
  \end{minipage}
\begin{figure}[H]
\caption{Symbolic eigenvalue computation times\protect\footnotemark$\;$for: Mathematica's built-in symbolic solver, and our 3 $\mathbb{Z}$LE libraries (Python, Mathematica, and C++).}
\label{fig:runtime}
\end{figure}
\footnotetext{Evaluated on a 2020 m1 MacBook Pro (16 GB, 8-core).}

\newpage
\section*{Applications}
$M^X$ matrices have elegant but constrained $\lambda \leftrightarrow a_{ij}$ structures. Arbitrary eigenvalue forms are reachable from $M^{X}$ via the spectral mapping theorem $f(\Lambda(A)) = \Lambda(f(A))$, tensors, wedge powers, and compounds\cite{2AC}, but most conveniently, direct modification $\lambda \rightarrow f(\lambda) \quad(a \rightarrow f(a))$. $\lambda = Ca$, so one can interchangeably work over $a$ or $\lambda$ for non-singular $C$. Given target $\mathbb{R}$-linear $\lambda'$: 
\[\lambda' = M\lambda\]
\[a' \rightarrow  C^{-1}\lambda' = C^{-1}MCa\]
Unlike forming random $Q^{-1}\Lambda'Q$, $A(a')$ inherits not only matrix features (symmetry, stochasticity, circulance, etc.) and parametric eigenvectors of $A$, but also its underlying algebraic structure $V=$span$\{B_i\}$.
\[A(a)= \begin{bmatrix}a_{1} & a_{2} & a_{3} \\a_{2} & a_{1} & a_{3} \\ a_{2} & a_{3} & a_{1}\end{bmatrix}=a_{1}\begin{bmatrix}1 & 0 & 0 \\ 0 & 1 & 0 \\ 0 & 0 & 1\end{bmatrix}+a_{2}\begin{bmatrix}0 & 1 & 0 \\ 1 & 0 & 0 \\ 1 & 0 & 0\end{bmatrix}+a_{3}\begin{bmatrix}0 & 0 & 1 \\ 0 & 0 & 1 \\ 0 & 1 & 0\end{bmatrix}= \sum_{j}^{}a_{j}B_{j}\]
$A(a')=\sum_{j}^{}a'_{j}B_{j} =\sum_{j}^{}a_{j}B'_{j}\;$. Let $P =C^{-1}MC$:\newline
\[A(Pa)=\sum_{j=1}^{n}Pa_{j}B_{j}=\sum_{i=1}^{n}\sum_{j=1}^{n}P_{ij}a_{j}B_{i} =\sum_{j=1}^{n}a_{j}\sum_{i=1}^{n}P_{ij}B_{i}=\sum_{j=1}^{n}a_{j}B'_{i} \quad (B'=BP^{T})\]
The linear reparametrization acts as a change of basis. Given the weak condition that $B_i$ are linearly independent it reduces to a map over $V$. That is, for $\Phi(a)\rightarrow \sum_{j}^{}a_{j}B_{j}$ ($\Phi:\mathbb{R}^n\rightarrow V, \quad \ker \Phi = \{0\}$), 
\[A(a') = T(A(a)),\quad  T=\Phi \circ C^{-1}MC \circ \Phi^{-1}.\]

\noindent Nonlinear forms via reparametrization of $M^{X}$ conserve the parametric eigenbasis but require stricter conditions to be treated as an operator transform $T(\Phi(a))$. To build multivariate polynomials in $\lambda$-space:
\[\lambda^T = \begin{bmatrix}\lambda_{1} & \lambda_{2} & \lambda_{3}\end{bmatrix}, \quad a^T = \begin{bmatrix}a_{1} & a_{2} & a_{3}\end{bmatrix}, \quad c^T = \begin{bmatrix}c_{1} & c_{2} & c_{3}\end{bmatrix},\quad M\in \mathbb{R}^{3x3}\]
\begin{minipage}{0.5\textwidth}
  \centering
\[\text{Let }f(\lambda, c)=e^{M\ln{(\lambda)}+ \ln{c}}\]
\[A(a\rightarrow C^{-1}f)\]
\[\lambda(A)_{i}= c_{i}\lambda_{1}^{m_{i1}}\lambda_{2}^{m_{i2}}\lambda_{3}^{m_{i3}}\]\newline
\end{minipage}
\begin{minipage}{0.5\textwidth}
  \centering
\[\text{Let }g=Lf(\lambda,c)\]
\[A(a\rightarrow C^{-1}g)\]
\[\lambda(A)_{i}= l_{ij}c_{j}\lambda_{1}^{m_{j1}}\lambda_{2}^{m_{j2}}\lambda_{3}^{m_{j3}}\]\newline
\end{minipage}
Matrices with highly non-trivial symbolic \hyperlink{wex}{structures} can be viewed as transformations over $M^{X}$, computed with assistance from the provided eigenvalue solver.\newpage

\noindent $\exists A$ of type \ref{fig:dbo} with $\lambda_A=Ca, \;\;x = \begin{bmatrix}x_{1} & x_{2} & \dots & x_{n}\end{bmatrix}, \;\; x_{i} \in \mathbb{R} \;\; s.t.$,
\[S =\begin{bmatrix}  \horzbar & \phi_{1}(x)  & \horzbar \\ \horzbar & \phi_{2}(x) & \horzbar\\ & \vdots &  \\ \horzbar & \phi_{n}(x) & \horzbar \end{bmatrix} = \begin{bmatrix} | & | & & | \\ \phi_{1}(x) & \phi_{2}(x) & \dots & \phi_{n}(x) \\ | & | & & | \end{bmatrix} \; (S_{ij} = S_{ji}),\] 
\[\lambda_{S A S^{T}} = (C^{-1}S^2C)\lambda_{A}\].
Equivalently, $\Lambda_{S A S^{T}} = f(\Lambda_{A})$ under the relabeling $f(a) : a \rightarrow aS^{2}.$\newline

\noindent The class of matrices whose eigenvalues can be computed using $\mathbb{Z}$-linear matrix computations as a subroutine is nascent and worth exploring. For example, the above construction\footnote{This construction does not guarantee $\mathbb{R}$-linear eigenvalues of the form $f(\Lambda_{A})$ for all matrices. See the \hyperlink{cons}{appendix} for examples and exceptions.}can yield  $\mathbb{R}$-linear matrices through a fractional transform\cite{FGT} on $A$ with bistochastic $S$, whose rows (columns) are permutations $\phi$ of some real values $x$.\newline

\noindent Furthermore, we can treat stochastic $\mathbb{Z}$-linear eigenvalue matrices as generator matrices via row-sum subtraction, and adopt a multidimensional Markov chain formalism\cite{MDMC} often used in stochastic automata networks (SAN).\newline

\begin{figure}[H]
\[ Q = Q_{0} + Q_{D}, \;\; Q_0 = \sum_{k=1}^{K}\bigotimes_{h = 1}^{H} Q_{k}^{(h)}, \;\; Q_{D} = -\text{diag}{(Q_{0}e)}\]
\[Q_{local} =\sum_{k=1}^{H}\bigotimes_{h = 1}^{H} Q_{k}^{(h)}, \; Q_{synchronized} = \sum_{k=H+1}^{K}\bigotimes_{h = 1}^{H} Q_{k}^{(h)}\]
\caption{MDMC. $e$: the $\mathbf{1}$s vector, $H$: $\#$ factors in each term, $K$: $\#$ terms. Subsystems $1, \hdots, H$ local transitions, $H+1, \hdots, K$ synchronized transitions.}
\label{fig:oooo}
\end{figure}
\noindent Such models are relevant because of the properties of Kronecker products and sums\cite{tp}. For eigenvalues $\lambda$ and $\mu$ of $A$ and $B$ and their eigenvectors $x$ and $y$, $A \otimes B$ has an eigensystem $\{\lambda \mu, \;x \otimes y\}$, and $A \oplus B$ $\{\lambda + \mu, \;x \otimes y\}$. Since $P=e^{Q}Q(0)$, we can see how the term 'local' comes 
from the capacity to factor the exponential and isolate the first $H$ terms, while synchronized transitions couple due to the eigenvalues being products rather than sums. By treating $M^{X}$ matrices as building blocks in the MDMC formalism we can easily trace spectral contributions from each factor.  Consider two matrices $A$ and $B$ - 
\newline\\
$A = \begin{bNiceMatrix}[first-row,last-row,first-col,last-col]
 & _{1} & _2 & _3 & \\
_1 & a & b & c & \\
_2 & b & a & c & \\ 
_3 & c & b & a & \\
 &   &   &   & \end{bNiceMatrix}\;\;\;$ \scalebox{0.51}{
 \begin{tikzpicture}
    \node[shape=circle,draw=black] (A) at (0,0) {1};
    \node[shape=circle,draw=black] (B) at (4,0) {2};
    \node[shape=circle,draw=black] (C) at (2,2*1.7321) {3};
        
    \path[<->,draw=black]
    (A) edge node[below] {b} (B)
    (C) edge node[left] {c} (A);

    \path[arrows = {->[harpoon]},transform canvas={xshift=-1.6,yshift=-1.6}]
    (B) edge node [at end,left, below,yshift=-4] {b} (C);
    \path[arrows = {->[harpoon]},transform canvas={xshift=1.6,yshift=1.6}]
    (C) edge node [at end,right, above, yshift=4] {c} (B);

    \path[->,draw=black]
    (A) edge [out=240,in=210,looseness=8] node[below] {a} (A)
    (B) edge [out=330,in=300,looseness=8] node[below] {a} (B)
    (C) edge [loop above] node {a} (C);
\end{tikzpicture}
}
$\qquad\;\; B = \begin{bNiceMatrix}[first-row,last-row,first-col,last-col]
 & _4 & _5 & \\
_4 & d & e & \\ 
_5 & e & d & \\
 &   &   & \end{bNiceMatrix}$
\scalebox{0.52}{
\begin{tikzpicture}
    \node[shape=circle,draw=black] (A) at (0,0) {1};
    \node[shape=circle,draw=black] (B) at (4,0) {2};
        
    \path[<->,draw=black]
    (A) edge node[below] {e} (B);
      
    \path[->,draw=black]
    (A) edge [out=75,in=105,looseness=8] node[above] {d} (A)
    (B) edge [out=75,in=105,looseness=8] node[above] {d} (B);
\end{tikzpicture}
}
\begin{minipage}{0.6\textwidth}

\[a+b+c+d+e=k\]
\[Q_{A \oplus B} = \begin{bNiceArray}{ccccccc}[small, first-col,first-row]
& _{(1,4)} & _{(1,5)} & _{(2,4)} & _{(2,5)} & _{(3,4)} & _{(3,5)} & \\
_{(1,4)} & a + d -k & e & b & 0 &c & 0 & \\ 
_{(1,5)} & e & a + d -k & 0 & b & 0 & c &\\
_{(2,4)} & b & 0 & a + d -k & e & c & 0 &\\
_{(2,5)} & 0 & b & e & a + d -k & 0 & c & \\
_{(3,4)} & c & 0 & b & 0 & a + d -k & e &\\
_{(3,5)} & 0 & c & 0 & b & e & a + d-k & \\
 & & & & & & \end{bNiceArray}\]
\[\sigma_{Q_{A \oplus B}}=-k + \text{map}((\lambda,\mu) \rightarrow \lambda + \mu,\;\sigma_A \bigtimes \sigma_B)\]
\end{minipage}
\begin{minipage}{0.4\textwidth}\raggedright
    \scalebox{0.51}{
    \begin{tikzpicture}
    \node[shape=circle,draw=black] (A) at (0,0) {(1,4)};
    \node[shape=circle,draw=black] (B) at (4,0) {(2,4)};
    \node[shape=circle,draw=black] (C) at (2,2*1.7321) {(3,4)};

    \node[shape=circle,draw=black] (A2) at (3*1.73,3) {(1,5)};
    \node[shape=circle,draw=black] (B2) at (4+3*1.73,3) {(2,5)};
    \node[shape=circle,draw=black] (C2) at (2+3*1.73,2*1.7321+3) {(3,5)};

    \path[<->,draw=black]
    (A) edge node[below] {e} (A2)
    (B) edge node[below] {e} (B2)
    (C) edge node[below] {e} (C2)
    (A) edge node[below] {b} (B)
    (C) edge node[left] {c} (A)
    (A2) edge node[below] {b} (B2)
    (C2) edge node[left] {c} (A2);

    \path[arrows = {->[harpoon]},transform canvas={xshift=-1.6,yshift=-1.6}]
    (B) edge node [at end,left, below,yshift=-4] {b} (C)
    (B2) edge node [at end,left, below,yshift=-4] {b} (C2);
    \path[arrows = {->[harpoon]},transform canvas={xshift=1.6,yshift=1.6}]
    (C) edge node [at end,right, above, yshift=4] {c} (B)
    (C2) edge node [at end,right, above, yshift=4] {c} (B2);
        
    \path[->,draw=black]
    (A) edge [out=240,in=210,looseness=8] node[below] {a+d} (A)
    (B) edge [out=330,in=300,looseness=8] node[below] {a+d} (B)
    (C) edge [loop above] node {a+d} (C)
    (A2) edge [out=240,in=210,looseness=8] node[below] {a+d} (A2)
    (B2) edge [out=330,in=300,looseness=8] node[below] {a+d} (B2)
    (C2) edge [loop above] node {a+d} (C2);
    \end{tikzpicture}
    }
\end{minipage}
\begin{minipage}{0.6\textwidth}
\[(a+b+c)(d+e)=k\]
\[Q_{A \otimes B} = \begin{bNiceArray}{ccccccc}[small, first-col,first-row]
& _{(1,4)} & _{(1,5)} & _{(2,4)} & _{(2,5)} & _{(3,4)} & _{(3,5)} & \\
_{(1,4)} & a d -k & a e & b d & b e & c d & c e \\
_{(1,5)} & a e & a d -k & b e & b d & c e & c d \\
_{(2,4)} & b d & b e & a d -k & a e & c d & c e \\
_{(2,5)} & b e & b d & a e & a d -k & c e & c d \\ 
_{(3,4)} & c d & c e & b d & b e & a d -k & a e \\
_{(3,5)} & c e & c d & b e & b d & a e &a d -k \\
& & & & & & \end{bNiceArray}\]

\[\sigma_{Q_{A \otimes B}}=-k + \text{map}((\lambda,\mu) \rightarrow \lambda \cdot \mu,\;\sigma_A\bigtimes \sigma_B)\]
\end{minipage}
\begin{minipage}{0.4\textwidth}
    \scalebox{0.42}{
    \begin{tikzpicture}
    \node[shape=circle,draw=black] (A2) at (5*1,0) {(1,5)};
    \node[shape=circle,draw=black] (C2) at (5*0.5,5*0.866) {(3,5)};
    \node[shape=circle,draw=black] (C) at (-0.5*5,0.866*5) {(3,4)};
    \node[shape=circle,draw=black] (A) at (-1*5,0) {(1,4)};
    \node[shape=circle,draw=black] (B) at (-0.5*5,-0.866*5) {(2,4)};
    \node[shape=circle,draw=black] (B2) at (0.5*5,-0.866*5) {(2,5)};

    \path[<->,draw=black]
    (A) edge node[left] {bd} (B)%
    (A) edge node[left] {cd} (C)%
    (A) edge node[xshift=2.2em,yshift=0.3em] {ae} (A2)%
    (A) edge node[below,xshift=-2em,yshift=1.4em] {be} (B2)%
    (A) edge node[below,xshift=2em,yshift=1.2em] {ce} (C2)%
        
    (B) edge node[below,xshift=2em,yshift=1.2em] {be} (A2)%
    (B) edge node[below] {ae} (B2)%
    
    (C) edge node[below,xshift=-2em,yshift=1.12em] {ce} (A2)%
    (C) edge node[above] {ae} (C2)%

    (A2) edge node[right] {bd} (B2)%
    (A2) edge node[right] {cd} (C2);%

    \path[->,draw=black]
    (B) edge [bend left] node[above,xshift=0.55em,yshift=0.2em] {cd} (C)%
    (C) edge node[above,xshift=-0.58em,yshift=0.2em] {bd} (B)%
    (B) edge [bend right] node[xshift=0.45em, yshift=-0.3em] {ce} (C2)%
    (C2) edge node[xshift=-2em, yshift=-2em] {be} (B)%
    (B2) edge [bend left] node[xshift=-0.5em, yshift=-0.3em] {ce} (C)%
    (C) edge node[xshift=-2em, yshift=2em] {be} (B2)%
    (B2) edge [bend right] node[above,xshift=-0.55em,yshift=0.2em] {cd} (C2)%
    (C2) edge node[above,xshift=0.58em,yshift=0.2em] {be} (B2)%

    (A) edge [loop left] node[below] {ad} (A)
    (B) edge [out=240,in=210,looseness=8] node[below] {ad} (B)
    (C) edge [loop above] node {ad} (C)
    (A2) edge [loop right] node[below] {ad} (A2)
    (B2) edge [out=330,in=300,looseness=8] node[below] {ad} (B2)
    (C2) edge [loop above] node {ad} (C2);
    \end{tikzpicture}
    }
\end{minipage}
\begin{figure}[H]
\caption{Local $A \oplus B$ (above) and synchronous $A \otimes B$ transitions (below). Eigenvalues $\sigma$ are sums or products of those from $A$ and $B$. State diagrams for each transition type are provided beside each $Q$-matrix.}
\label{fig:trans}
\end{figure}


\noindent SAN are expressive and have many memory- and time-efficient representations.
\[Q_0=\bigoplus_{i=1}^{H} Q_{i}^{(i)}\]
\[e^{Qt}Q(0) = e^{Q}=e^{\bigoplus_{i=1}^{H} Q_{i}^{(i)} + Q_{D}} = e^{\bigoplus_{i=1}^{H}Q_{i}^{(i)}}e^{Q_{D}}=e^{Q_{D}}\bigotimes_{h=1}^{H}e^{Q_{i}^{(i)}} = e^{Q_{D}}\bigotimes_{h=1}^{H}P_{i}^{(i)}\]
Relating $P$ and $Q$ in a local MDMC by exploiting the identity $e^{x \oplus y} = e^{x} \otimes e^{y}$\cite{fax}

\newpage

\begin{figure}[H]
{\includegraphics[width=.5\textwidth]{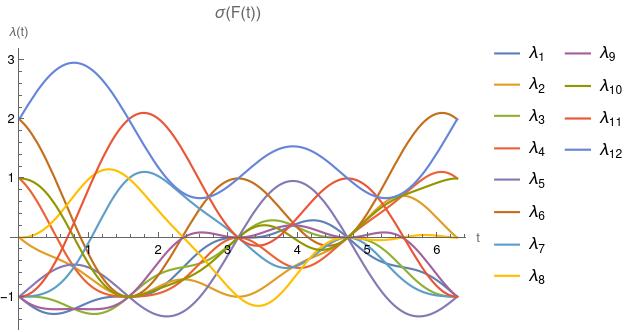}}
{\includegraphics[width=.49\textwidth]{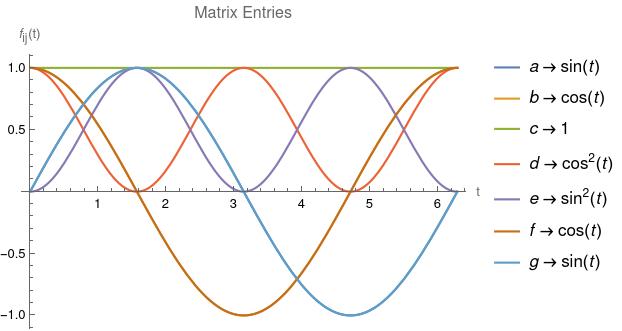}}
\caption{Given eigenvalue constraints for $F(t)$, find time-varying matrix entries $f_{ij}(t)$ ($\sigma(F(t))$ plotted with scale parameter $s=\frac{1}{2}$).}

\label{fig:ooooo}
\end{figure}
\[F(t) = s\bigoplus_{i=1}^{H} Q_{i} +  (1-s)\bigotimes_{i=1}^{H} Q_{i} + Q_D\]
\[Q_1 = \begin{bsmallmatrix}a(t) & b(t) \\ b(t) & a(t)\end{bsmallmatrix}, \; Q_2 = \begin{bsmallmatrix}
c(t) & d(t) & e(t) \\
d(t) & c(t) & e(t) \\ 
e(t) & d(t) & c(t)
\end{bsmallmatrix}, \; Q_3 = \begin{bsmallmatrix}f(t) & g(t) \\ g(t) & f(t)\end{bsmallmatrix}\]
\newline
One toy model, $F(t)$, is a $12$-dimensional convex combination of local and synchronous transitions sharing the same 3 factors $Q_i$ of dimensions $2,3,2$. $s \in [0,1]$ scales their proportions. Note that in general the system is underdetermined but particular solutions exist, such as those for the eigenvalues shown in fig. \ref{fig:ooooo}.

\section*{Conclusion}
Users can generate matrices and find symbolic eigenvalues which are $\mathbb{Z}$-linear combinations of matrix entries. Specifically, users can generate $M^X$ with any $X$, alongside restricted stochastic $M^X$ with two distinct matrix generation scripts. They can compute fast symbolic eigenvalues using the language of choice for high dimension on a commercial laptop. Real-coefficient, multivariate eigenvalues can be achieved directly from $M^X$ calculations whilst preserving some structure. We hope this toolkit benefits researchers in any befitting domain, which includes but is not limited to reduced order modeling  \& control, conservative systems, and large-scale numeric simulations.

\section*{Future Directions}
A similar treatment for eigenvectors is equally important. Generating a closest suitable matrix still requires expert knowledge for distance norm selection and usage, so a procedural method for inferring best-fit $\mathbb{Z}$-linear matrices would be useful. Developed understanding of generalized eigenvalue structures is needed.  We suspect that more efficient encodings could decrease the required number of batches.  We look forward to progress in computing fast eigenvalues with complex algebraic structure through novel means of encoding symbolic terms in a numeric context.

\bibliography{cites}

\appendix
\section*{Appendix}
\hypertarget{alg}{\subsection*{Algorithm Details}}
\color{red}*\color{black} : undefined, pre-existing

\subsubsection*{Solver}
\vspace{-0.5cm} 
\begin{algorithm}[H]
\renewcommand\thealgorithm{}
\caption{$\lambda$ solver}
\begin{algorithmic}
\STATE $\text{A, symbols, batchsize, base, stagger} = \hdots$
\STATE $\text{precision, extractionparams (midpoint, midpointvalue), digits} = \hdots$
\STATE $.\text{indicatorvars} = \begin{bmatrix}\text{base}^{\ceil{\frac{\text{digits}}{2}} - \text{stagger}} & \text{base}^{\ceil{\frac{\text{digits}}{2}} - 2*\text{stagger}} & \hdots & \text{base}^{-\floor{\frac{\text{digits}} {2}}}\end{bmatrix}$
\STATE $.\text{realvals} = \begin{bmatrix} \hdots \end{bmatrix}$
\STATE $\text{map} = \text{\{symbols} \rightarrow \text{randomreals\}}$
\STATE $\text{coeffs} = \text{{[ ]}}$

\FOR{$i=0$ to \text{numbatches}}
\STATE $\text{mapcurr} = \text{map.copy()}$
\STATE $\text{mapcurr{[}}i*\text{batchsize}, \min{\text{(n}, (i+1)*\text{batchsize)}}\text{{]}} \mathrel{+}= \text{expvars}$
\STATE $\text{A}_{i} = \text{mapcurr(A)}$
\STATE $\sigma = \text{eig\color{red}*\color{black}}(\text{A}_{i}) + \text{midpoint}$
\STATE $\text{coeffs.insert({[} digitsof}(\sigma_{j}\text{)} - \text{midpointvalue for j in range(n){]})}$
\ENDFOR
\STATE $\text{coeffs} = \begin{bmatrix} \text{vec(coeffs[:},0\text{{]})} \\ \hdots \\ \text{vec(coeffs[:},n-1\text{{]})} \end{bmatrix}$

\end{algorithmic}
\end{algorithm}
\subsubsection*{General Matrix Generation}
\label{sec:gen}
\vspace{-0.5cm} 
\begin{algorithm}[H]
\renewcommand\thealgorithm{}
\caption{matrix-generation(partialorder)}
\begin{algorithmic}
\STATE $\phi = \text{get-phi(dcs\color{red}*\color{black}(graph(partialorder)))}$ 

\FOR{$i=0$ to \text{len}$(\phi)$}
\FOR{$j=1$ to \text{len}$(\phi)$}
\STATE $A[i,j] = \epsilon$\color{red}*\color{black}($\phi_{j}\phi_{i}^{-1}$) 
\ENDFOR
\ENDFOR
\STATE $\text{adstrings} = \text{uniques}(A)$
\STATE $\text{s} = [a_1, \dots, a_{\text{len(adstrings)}}]$

\RETURN $A.\{\text{adstrings} : \text{s}\}$
\end{algorithmic}
\end{algorithm}
\vspace{-0.5cm} 
\begin{algorithm}[H]
\renewcommand\thealgorithm{}
\caption{get-phi(dcs)}
\begin{algorithmic}
\STATE $\text{g} = [\text{map(transitive-reduction\color{red}*\color{black}, i) for i in dcs}]$
\STATE $\text{grouppos, groupnodes} = \text{flatten(}[\text{i.edges for i in g}]\text{), } [\text{i.nodes for i in g}]$
\STATE $\text{perms = fixed-perms(g, grouppos, groupnodes)}$ 
\STATE $\text{orderings = combination-labels(groupnodes)}$ 
\STATE $\phi = \text{perms} \varprod_{}^{} \text{orderings}$

\RETURN $\phi$
\end{algorithmic}
\end{algorithm}
\vspace{-0.5cm} 
\begin{algorithm}[H]
\renewcommand\thealgorithm{}
\caption{combination-labels(nodes, groupnodes)}
\begin{algorithmic}
\STATE $\text{nodes} = \bigcup\text{groupnodes}$
\STATE $\text{prod = combinations\color{red}*\color{black}(nodes, len(groupnodes}[0]\text{)}$
\FOR{$i=1$ to \text{len(groupnodes)}}
\STATE $\text{subsets = combinations\color{red}*\color{black}(nodes, len(groupnodes}[i]\text{)}$
\STATE $\text{matches} = \text{{[ ]}}$
\FOR{$j=0$ to \text{len(prod)}}
\FOR{$k$ in \text{subsets}}
\IF{$\text{len(k}\setminus\text{prod}[j]\text{) == len(groupnodes}[i]\text{)}$}
\STATE $\text{matches.insert(flatten(prev}[j]\text{, k} \setminus \text{prod}[j]\text{))}$
\ENDIF
\ENDFOR
\ENDFOR
\STATE $\text{factors} = \text{zip(prod, matches)}$
\STATE $\text{prod} = \text{map(}\varprod_{}^{}\text{, factors)}$
\ENDFOR
\STATE $\text{remainder} = [\text{nodes}\setminus\text{i for i in prod}]$
\RETURN $[\text{prod}[i]+\text{remainder}[i] \text{ for i in range(len(prod))}]$
\end{algorithmic}
\end{algorithm}

\vspace{-0.5cm} 
\begin{algorithm}[H]
\renewcommand\thealgorithm{}
\caption{fixed-perms(g, grouppos, groupnodes)}
\begin{algorithmic}
\STATE $\text{dc-perms} = \text{{[ ]}}$
\FOR{$i=0$ to \text{len(g)}}
\STATE $\text{nodesi} = \text{groupnodes}[i]$
\STATE $\text{n} = \text{len(nodesi)}$
\STATE $\text{nodemap} = \text{dict(zip(nodesi, range(len(nodesi))))}$
\STATE $\text{poi} = \text{grouppos}[i]$
\STATE $\text{abovebelow} = \text{bi-level(g}[i]\text{, nodesi)}$ 
\STATE $\text{permsi} = \text{get-perms(n, nodesi, abovebelow, levels, nodemap)}$ 
\STATE $\text{permsi} = \text{filter(permsi, i, nodesi, poi, nodemap)}$ 
\STATE $\text{dc-perms.insert(permsi)}$
\ENDFOR

\RETURN $\varprod_{}^{}\text{dc-perms}$
\end{algorithmic}
\end{algorithm}
\vspace{-0.5cm} 
\begin{algorithm}[H] 
\renewcommand\thealgorithm{}
\caption{bi-level(g, nodes)}
\begin{algorithmic}
\STATE $\text{rg} = \text{reverse-graph\color{red}*\color{black}(g)}$
\STATE $\text{tg} = \text{topological-sort\color{red}*\color{black}(g)}$
\STATE $\text{trg} = \text{topological-sort\color{red}*\color{black}(rg)}$
\STATE $\text{abovebelow} = \text{{[ ]}}$

\FOR{$\text{node}$ in \text{nodes}}
\STATE $\text{above} = \text{strongconnect\color{red}*\color{black}([\text{node}], \text{g}, \text{tg})}$
\STATE $\text{below} = \text{strongconnect\color{red}*\color{black}([\text{node}], \text{rg}, \text{tgr})}$
\STATE $\text{abovebelow}[\text{node}] = [\text{len(above)}-1, \text{len(below)}-1]$
\ENDFOR
\RETURN $\text{abovebelow}$
\end{algorithmic}
\end{algorithm}
\vspace{-0.5cm} 
\begin{algorithm}[H]
\renewcommand\thealgorithm{}
\caption{get-perms(n, nodes, abovebelow, levels, nodemap)}
\begin{algorithmic}
\STATE $\text{legalvalues} = [\text{nodes}[0+i[0]:n-i[1]] \text{ for i in abovebelow.values}()] $
\STATE $\text{blockvariables}=[[\text{map(legalvalues[}\#\text{.node-map], i)}] \text{ for i in levels.values()}]$
\STATE $\text{blocks} = \text{{[ ]}}$
\FOR{$i$ in \text{blockvariables}}
\STATE $\text{blocks.append(filter-duplicates}(\varprod_{}^{}i))$
\ENDFOR
\STATE $\text{res} = \text{filter-duplicates\color{red}*\color{black}}(\varprod_{}^{}\text{blocks})$
\STATE $\text{res} = \text{[i for i in res if len(uniques(i))=len(i)]}$
\STATE $\text{reorder} = \text{perm-inverse\color{red}*\color{black}(}\{\text{nodemap.values(), flatten(levels.values())}\}\text{)}$
\RETURN $[i[\text{reorder}] \text{ for i in res}]$
\end{algorithmic}
\end{algorithm}
\vspace{-0.5cm} 
\begin{algorithm}[H]
\renewcommand\thealgorithm{}
\caption{filter(perms, i, nodes, po, nodemap)} 
\begin{algorithmic}
\STATE $\text{indices} = \text{{[ ]}}$
\FOR{$i=0$ in \text{perm in perms}}
\FOR{$j=0$ in $\text{range(}0,\text{int(}\frac{\text{len(po)}}{2}\text{))}$}
\STATE $\text{relation} = \text{order}[2j:2(j+1)]$
 \IF{$\text{perms}[i][\text{nodemap}[\text{relation}[0]]] < \text{perms}[i][\text{nodemap}[\text{relation}[1]]]$}
     \STATE $\text{indices.insert(}i\text{)}$
\ENDIF
\ENDFOR
\ENDFOR
\RETURN $\text{perms}[\text{indices}]$
\end{algorithmic}
\end{algorithm}
\vspace{-0.5cm} 


\subsubsection*{Disjoint, Local Transposition Matrix Generation}
\label{sec:fgen}
\vspace{-0.5cm} 
\begin{algorithm}[H]
\renewcommand\thealgorithm{}
\caption{fib-matrix-generation(fac)}
\begin{algorithmic}
\STATE $\text{adstrings} = \text{ascent-descent(fac)}$ 
\FOR{$i=0$ to \text{len(adstrings)}}
\FOR{$j=0$ to \text{len(adstrings)}}
\STATE $A[i,j] = \text{ad-mult(adstrings[i], adstrings[j])}$ 
\ENDFOR
\ENDFOR

\RETURN $A.\{\text{adstrings} : [x_{1}, \dots, x_{\text{len(adstrings)}}]\}$
\end{algorithmic}
\end{algorithm}
\vspace{-0.5cm} 
\begin{algorithm}[H]
\renewcommand\thealgorithm{}
\caption{ascent-descent(fac)}
\begin{algorithmic}[]
\STATE $\text{segments} = \text{map(reduce-factors, fac)}$ 
\STATE $\text{map(append-0\color{red}*\color{black}, segments)}[0:len(fac)-1]$ 
\STATE $\text{res} = \varprod_{i}^{}\text{segments}[i]$ 

\RETURN $\text{res}$
\end{algorithmic}
\end{algorithm}
\vspace{-0.5cm} 

\begin{algorithm}[H]
\renewcommand\thealgorithm{}
\caption{ad-mult(s1,s2)}
\begin{algorithmic}[]
\STATE $\text{res} = \text{s1.copy()}$
\STATE $\text{z1, z2} = \text{where(s1 = 0), where(s2 = 0)}$
\STATE $\text{difs} = [\text{abs}(i[0] - i[1]) \text{ for i in }  \text{z1}\varprod_{}^{}\text{z2}]$ 
\STATE $\text{zdif, odif} = \text{where(difs = 0), where(difs = 1)}$
\STATE $\text{marks} = \text{zdif} \cup \text{odif}$
\FOR{$i$ in \text{marks}}
\STATE $\text{res}[i] = 1$
\ENDFOR
\FOR{$i$ in \text{z2 }$\setminus$\text{ zdif}}
\STATE $\text{res}[i] = 0$
\ENDFOR

\RETURN \text{res}
\end{algorithmic}
\end{algorithm}
\vspace{-0.5cm} 
\begin{algorithm}[H]
\renewcommand\thealgorithm{}
\caption{reduce-factors(x, lists)} 
\begin{algorithmic}[]

\STATE $\text{lists} = [[1],[0]]$
\FOR{$i=0$ to \text{x}} 
\STATE $\text{sub} = \text{{[ ]}}$
\FOR{$j$ in \text{lists}}

 \IF{$j[-1] == 1$}
     \STATE $\text{sub.insert}(j+[0])$
\ENDIF

\STATE $\text{sub.insert}(j+[1])$

\ENDFOR
\STATE $\text{lists} = \text{sub}$
\ENDFOR

\RETURN $\text{lists}$
\end{algorithmic}
\end{algorithm}
\newpage

\hypertarget{cons}{\subsection*{Construction Examples}}
\noindent $\bold{n=3}$
\[S=\begin{bmatrix}x_{1} & x_{2} & x_{3} \\ 
x_{2} & x_{3} & x_{1} \\
x_{3} & x_{1} & x_{2}\end{bmatrix}\; A = \begin{bmatrix} a_{1} & a_{2} & a_{3} \\ a_{2} & a_{1} & a_{3} \\ a_{3} & a_{2} & a_{1} \end{bmatrix}, \;\sigma_A=\begin{bmatrix}1 & -1 & 0 \\ 1 & 0 & -1 \\ 1 & 1 & 1\end{bmatrix}\begin{bmatrix}a_{1} \\ a_{2} \\ a_{3}\end{bmatrix}\]
\[\sigma_{SAS^{T}} =C^{-1}S^2C\begin{bmatrix} a_{1} - a_{2} \\ a_{1} - a_{3} \\ a_{1} + a_{2} + a_{3}\end{bmatrix}, \quad C^{-1}S^2C=\text{diag}(\begin{bmatrix} x_{1}^2 - x_{1}x_{2} + x_{2}^2 - x_{1}x_{3} - x_{2}x_{3} + x_{3}^2 \\ x_{1}^2 - x_{1}x_{2} + x_{2}^2 - x_{1}x_{3} - x_{2}x_{3} + x_{3}^2 \\ (x_{1} + x_{2} + x_{3})^{2}\end{bmatrix})\]\newline
$\bold{n=5}$ 
\[S = \begin{bmatrix}\horzbar \phi_{1}(
x_{0} & x_{1} & x_{1} & x_{1} & x_{1}) \horzbar \\ & &\vdots & & \\\horzbar \phi_{5}(
x_{0} & x_{1} & x_{1} & x_{1} & x_{1}) \horzbar\end{bmatrix}, \quad A = \begin{bmatrix}a_{1} & a_{2} & a_{3} & a_{4} & a_{5} \\
a_{2} & a_{1} & a_{3} & a_{5} & a_{4}\\ 
a_{3} & a_{2} & a_{1} & a_{4} & a_{5} \\
a_{4} & a_{5} & a_{3} & a_{1} & a_{2} \\ 
a_{5} & a_{4} & a_{3} & a_{2} & a_{1}\end{bmatrix}\]
Requires added constraints on $x$, such as $\sum_{}^{}g(A_{*i}) = k$ with $g(a) : a \rightarrow \phi(x)$.\newline

\noindent $\bold{n=6}$
\[S=\begin{bmatrix}x_{1} & x_{2} & x_{3} & x_{4} & x_{5} & x_{6}\\  
x_{2} & x_{1} & x_{4} & x_{3} & x_{6} & x_{5} \\
x_{3} & x_{4} & x_{5} & x_{6} & x_{1} & x_{2} \\
x_{4} & x_{3} & x_{6} & x_{5} & x_{2} & x_{1} \\
x_{5} & x_{6} & x_{1} & x_{2} & x_{3} & x_{4} \\
x_{6} & x_{5} & x_{2} & x_{1} & x_{4} & x_{3} \end{bmatrix}, 
\;\; 
A=\begin{bmatrix}a_{1} & a_{2} & a_{3} & a_{4} & a_{5} & a_{6} \\
a_{2} & a_{1} & a_{4} & a_{3} & a_{6} & a_{5} \\
a_{3} & a_{4} & a_{1} & a_{2} & a_{5} & a_{6} \\ 
a_{4} & a_{3} & a_{2} & a_{1} & a_{6} & a_{5} \\ 
a_{5} & a_{6} & a_{3} & a_{4} & a_{1} & a_{2} \\ 
a_{6} & a_{5} & a_{4} & a_{3} & a_{2} & a_{1} \end{bmatrix}\]
\[=\begin{bmatrix} a_{1} & a_{2} & a_{3} \\ a_{2} & a_{1} & a_{3} \\ a_{3} & a_{2} & a_{1} \end{bmatrix}\otimes\begin{bmatrix}j & k \\ k & j\end{bmatrix}\]
$\sigma_A=[a_{1} + a_{2} - a_{3} - a_{4}, \; a_{1} - a_{2} - a_{3} + a_{4}, \;a_{1} + a_{2} - a_{5} - a_{6}, \; a_{1} - a_{2} + a_{3} - a_{4} + a_{5} - a_{6}, \;a_{1} - a_{2} - a_{5} + a_{6},\; a_{1} + a_{2} + a_{3} + a_{4} + a_{5} + a_{6}]$
\[\text{e.g. } x=[0.4\; 0.5\; -1\; 0.3\; .2\;.6], \; \sigma_{SAS^{T}} = [1.17\Lambda_{A_{1}}\; 1.17\Lambda_{A_{2}}\;2.41\Lambda_{A_{3}}\;2.41\Lambda_{A_{4}}\; 3.24\Lambda_{A_{5}}\; \Lambda_{A_{6}}]\]\newline
\noindent $\bold{n=6}$
\[S=\begin{bmatrix}x_{1}^{3x3} & x_{2}^{3x3} \\ x_{2}^{3x3} & \pm x_{1}^{3x3}\end{bmatrix}, \;\;A = \begin{bmatrix}B & D \\ D& B\end{bmatrix}, \;\; B,D\:=A(a)^{3x3}, A(a')^{3x3}\]
\[\frac{1}{2}\begin{bmatrix}1 & 1 \\ 1 & -1\end{bmatrix}\begin{bmatrix}B & D \\ D & B\end{bmatrix}\begin{bmatrix}1 & 1 \\ 1 & -1\end{bmatrix}^{T} = \begin{bmatrix}B+D & 0 \\ 0 & B-D\end{bmatrix}\]
\[\sigma_{A} = \sigma_{B+D}\cup\sigma_{B-D}\]
\[\sigma_{S A S^{T}} = (C^{-1}S^2C)\sigma_{A}\]

\hypertarget{wex}{\subsection*{Generic Eigenvalues}}
\noindent Consider the target symbolic eigenvalues 
\[\begin{bmatrix}\frac{a_{1} - a_2 - a_3}{2a_{2}} &  a_2 + \frac{a_{1} - a_2 - a_3}{a_{2}}& 3a_{2}a_3 + a_2\end{bmatrix}\]
One function composition route: \newline\\
$(1)$ identify linear terms via $M\lambda$, \newline\\
$(2)$ build multivariate mixtures via $f(M\lambda)$, and \newline\\
$(3)$ take sums over mixed terms via $g(f(M\lambda))$. \newline
\[A= \begin{bmatrix}a_{1} & a_{2} & a_{3} \\a_{2} & a_{1} & a_{3} \\ a_{2} & a_{3} & a_{1}\end{bmatrix}, \;\; \lambda^T =[a_{1}-a_{2},\; a_{1}-a_{3}, \; a_{1}+a_{2}+a_{3}]\]\newline\\
$(1)$ Linear terms: $\begin{bmatrix}a_{1} - a_2 - a_3 &  a_2 & a_3\end{bmatrix}$. Map $\lambda_{A}$ from $A^{3x3}\in M^{X}$:
\[\begin{bmatrix} m_{11} & m_{12} & m_{13} \\ m_{21} & m_{22} & m_{23} \\ m_{31} & m_{32} & m_{33}\end{bmatrix}\begin{bmatrix}\lambda_{1} \\ \lambda_{2} \\ \lambda_{3}\end{bmatrix} = \begin{bmatrix}a_{1} - a_2 - a_3\\ a_2\\ a_3\end{bmatrix}\] 
\[\begin{matrix}m_{13} = \frac{-m_{11}(a_{1}-a_{2})- m_{12}(a_{1}-a_{3}) + (a_{1}-a_{2}-a_{3})}{a_{1}+a_{2}+a_{3}} \\ m_{22} = \frac{-m_{21}(a_{1}-a_{2}) + a_{2} -m_{23}(a_{1}+a_{2}+a_{3})}{a_{1}-a_{3}} \\ m_{32} =-\frac{m_{31}(a_{1}-a_{2}) + a_{3} - m_{33}(a_{1}+a_{2}+a_{3})}{a_{1}-a_{3}}\end{matrix}\]
with 6 free $m_{ij}$ variables. Let $\lambda' = M\lambda$:\newline\\ \newline
$(2)$ $f(\lambda',c)=e^{M'\ln{(\lambda')}+ \ln{c}}$ with $M'=\begin{bmatrix}1 & -1 & 0 \\ 0 & 1 & 0 \\ 0 & 1 & 1\end{bmatrix}$, $c=\begin{bmatrix} \frac{1}{2} \\ 1 \\ 3\end{bmatrix}$\newline
$ = \begin{bmatrix}\frac{\lambda'_{1}}{2\lambda'_{2}} & \lambda'_{2} & 3\lambda'_{2}\lambda'_{3}\end{bmatrix} =  \begin{bmatrix}\frac{a_{1} - a_2 - a_3}{2a_{2}} &  a_2 & 3a_{2}a_3\end{bmatrix}$\newline\\ \newline
$(3)$ Apply $g=Lf$ with $L=\begin{bmatrix}1 & 0 & 0 \\ 2 & 1 & 0 \\ 0 & 1 & 1\end{bmatrix}$ for target values $g(L \circ f)(M\lambda_A):$\newline
\[\sigma_{A(C^{-1}g)}=\begin{bmatrix}\frac{a_{1} - a_2 - a_3}{2a_{2}} &  a_2 + \frac{a_{1} - a_2 - a_3}{a_{2}}& 3a_{2}a_3 + a_2\end{bmatrix}\].



\end{document}